\documentclass[11pt,reqno,oneside]{amsart}

\usepackage{graphicx}    
\usepackage{amsfonts}
\usepackage{verbatim}
\usepackage{amssymb}
\usepackage{amsmath}
\usepackage[all]{xy}
\usepackage{hyperref}
\usepackage[numbers,sort&compress]{natbib}

\usepackage[it,small]{caption}  

\newcommand{\N}{{\mathbb N}}

\newcommand{\R}{{\mathbb R}}

\newcommand{\E}{\mathcal{E}}
\DeclareMathOperator*{\argmin}{arg\,min}

\newcommand{\residual}{{\mathcal R}}
\newcommand{\half}{{\textstyle \frac{1}{2}}}
\newcommand{\fourth}{{\textstyle \frac{1}{4}}}

\renewcommand{\vec}[1]{\mathbf{#1}}

\newtheorem{lemma}{Lemma}[section]

\numberwithin{equation}{section}

\allowdisplaybreaks[4]

\begin{document}

\title[Adaptive Refinement for the Quasicontinuum Approximation]
{Goal-Oriented Adaptive Mesh Refinement
for the Quasicontinuum Approximation
of a Frenkel-Kontorova Model}
\author{Marcel Arndt}
\author{Mitchell Luskin}
\begin{abstract}
The quasicontinuum approximation \cite{TadmorOrtizPhillips:1996} is a method to reduce
the atomistic degrees of freedom of a crystalline solid
by piecewise linear interpolation from
representative atoms that are nodes for a finite
element triangulation.
In regions of the crystal with a highly nonuniform deformation such
as around defects, every atom must be
a representative atom to obtain sufficient accuracy, but
the mesh can be coarsened away from such regions to remove
atomistic degrees of freedom while retaining sufficient
accuracy.  We present an error estimator and a related
adaptive mesh refinement algorithm for the quasicontinuum
approximation of a generalized Frenkel-Kontorova model that
enables a quantity of interest to be efficiently computed
to a predetermined accuracy.
\end{abstract}

\keywords{Error estimation, {\em a posteriori}, adaptive,
refinement, goal-oriented, atomistic-continuum modeling,
quasicontinuum, Frenkel-Kontorova model, dislocation, defect}

\subjclass[2000]{65Z05, 70C20, 70G75}

\thanks{This work was supported in part by DMS-0304326 and by the Minnesota
  Supercomputing Institute. This work is also based on work supported by the
  Department of Energy under Award Number DE-FG02-05ER25706.}


\maketitle


\section{Introduction}
The solution of the equations for mechanical equilibria of a crystalline solid
modeled by a classical atomistic potential requires the computation of the
interaction of each atom with all of the other atoms in its sphere of
influence.  Due to the high computational complexity, it is generally not
possible to obtain numerical solutions for systems that are large enough to
simulate long-range elastic effects, even for short-ranged potentials.  However,
the local environment of nearby atoms is almost identical up to translation,
except in the neighborhood of defects such as cracks and dislocations.  The
quasicontinuum method utilizes this slow variation of the strain away from
defects to approximate the full systems of equations of mechanical equilibrium
by equations of equilibrium at a reduced set of representative atoms
~\cite{MillerTadmor:2002,CurtinMiller:2003,TadmorOrtizPhillips:1996,DobsonLuskin:2006}.

More precisely, the positions of the full set of atoms are obtained by piecewise
linear interpolation from representative atoms that are nodes for a finite
element triangulation.  A quasicontinuum energy is defined as a function of the
positions of the representative atoms.  The quasicontinuum method is then used
to obtain a solution to a desired accuracy with a significant reduction in the
computational degrees of freedom by coarsening the finite element mesh away from
the singularities.  Near defects, sufficient accuracy can only be obtained if
all atoms are representative atoms.  In contrast to continuum models, the mesh
cannot be refined past the atomistic scale.

Even higher efficiency is achieved by approximating the total energy of the
atoms in a coarsened triangle by the product of the strain energy density and
the area of the triangle (or its higher dimensional analogue).  The strain energy
density is obtained from the energy per atom in a lattice which is a uniform
strain of the infinite reference lattice.  The uniform strain in turn
is determined from the displacement of the nodes of the respective triangle.
The quasicontinuum approximation we obtain this way allows the coupling between
a region that is computed as in fully atomistic simulations and a region that is
computed using the methods of piecewise linear finite element continuum
mechanics.

It is well-known that an efficient adaptive algorithm is highly dependent on the
quantity of interest or goal of the computation, see \cite{AinsworthOden:2000,
  BangerthRannacher:2003}.  In this paper, we adopt this goal-oriented approach
to obtain the approximation of a quantity of interest to within a desired
tolerance.

The reliable and efficient utilization of the quasicontinuum method requires
both a strategy to determine the decomposition of the fully atomistic system
into atomistic and continuum regions and a strategy for refinement within the
continuum region.  In \cite{ArndtLuskin:2007a, ArndtLuskin:2007b}, we have
developed a goal-oriented error estimator and a corresponding adaptive algorithm
to decide between the atomistic model and the continuum model.  In this paper,
we develop an error estimator and a corresponding adaptive algorithm for mesh
refinement in the continuum region.

We have chosen initially to investigate as a model problem a generalization of
the classical Frenkel-Kontorova model.  The potential energy includes
next-nearest-neighbor interactions in addition to nearest-neighbor interactions
so that the continuum energy of a representative atom is different from the
atomistic energy of a representative atom for a nonuniform strain.

We note that in contrast to error estimators and adaptive algorithms for mesh
refinement in classical continuum mechanics, our continuum region is coupled to
the atomistic region. This is a considerably more complex setting than in purely
continuum models with classical boundary conditions.  Also, our mesh refinement
algorithm restricts the representative atoms which serve as the mesh points to
the sites of the atoms in the reference lattice, although this could be relaxed
away from the atomistic-continuum interfacial region.

Algorithms for adaptive mesh refinement for the quasicontinuum method have been
proposed and investigated in numerical experiments for several mechanics
problems in \cite{MillerTadmorPhillipsOrtiz:1998,
  ShenoyMillerTadmorRodneyPhillipsOrtiz:1999, KnapOrtiz:2001}.  An {\em a
  posteriori} error indicator for a global norm was analyzed and tested for a
variant of our one-dimensional quasicontinuum method in
\cite{OrtnerSueli:2006b}.  The goal-oriented approach to adaptive mesh
refinement for the quasicontinuum method was first investigated in
\cite{OdenPrudhommeBauman:2006,OdenPrudhommeBauman:2005}.  Mathematical analyses
of several variants of the quasicontinuum method have been given in
\cite{BlancLeBrisLegoll:2005, ELuYang:2006, EMing:2004, Lin:2003, OrtnerSueli:2006a, Lin:2005}.  We refer to
\cite{BadiaParksBochevGunzburgerLehoucq:2007, ParksBochevLehoucq:2007,
  LiuKarpovZhangPark:2004, XiaoBelytschko:2004, CurtinMiller:2003} for
alternative atomistic-continuum coupling methods.


\section{Quasicontinuum Approximation}

In this section, we introduce the atomistic model and its
quasicontinuum approximation.  We refer to \cite{ArndtLuskin:2007a}
for a more detailed description.

We consider a one-dimensional system of $2M$ atoms whose positions are denoted
by $\vec y= (y_{-M+1},\dots, y_M)\in\R^{2M}$. The potential energy of the
atomistic system is described by a function
\begin{align}
  \E^a : \R^{2M} \to \R.
\end{align}
We split the energy into atom-wise contributions by means of
\begin{align} \label{EqADecomp}
  \E^a(\vec y) = \sum_{i=-M+1}^M \E^a_i(\vec y).
\end{align}

\begin{figure}
\centering
\includegraphics[width=0.8\textwidth]{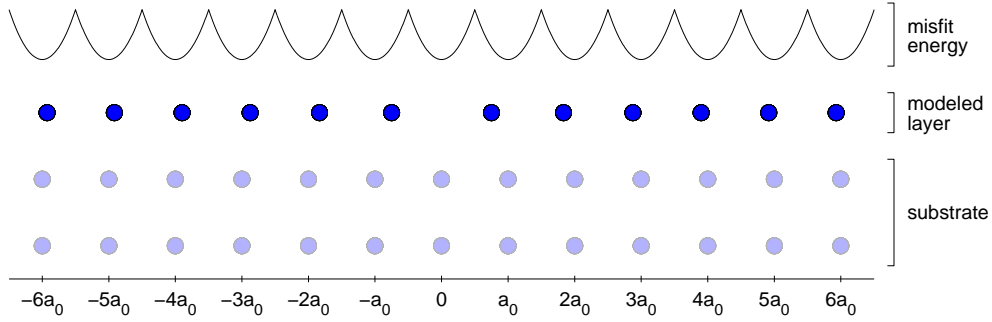}
\caption{Frenkel-Kontorova model. The wells depict the misfit energy,
  $\E^{m,a}_i$.}
\label{FigFKModel}
\end{figure}

In this paper, we consider a Frenkel-Kontorova type model \cite{Marder:2000}
which serves as a one-dimensional description of a dislocation.  We expect that
the {\em a posteriori} error estimators that we introduce and the corresponding
adaptive refinement algorithms will be applicable to more general quasicontinuum
models.  For the Frenkel-Kontorova model, the atom-wise contribution, $\E^a_i$,
consists of two parts,
\begin{align}
  \E^a_i = \E^{a,e}_i + \E^{a,m}_i.
\end{align}
The elastic part, $\E^{a,e}_i$, describes nearest-neighbor (NN) interactions and
next-nearest-neighbor (NNN) interactions, whereas $\E^{a,m}_i$ models the misfit
energy of a slip plane sitting on an undeformed substrate, see
Figure~\ref{FigFKModel}.  The two parts are defined as
\begin{equation} \label{EqDefEai}
\begin{split}
  \E^{a,e}_i(\vec y) & =
              \fourth k_1(y_i     - y_{i-1} -  a_0)^2
            + \fourth k_1(y_{i+1} - y_i     -  a_0)^2 \\
    & \quad + \fourth k_2(y_i     - y_{i-2} - 2a_0)^2
            + \fourth k_2(y_{i+2} - y_i     - 2a_0)^2,\\
  \E^{a,m}_i(\vec y) & = \begin{cases}
    \half k_0 \left( y_i - (i-1) a_0 \right)^2,  & i=-M+1,\ldots,0, \\
    \half k_0 \left( y_i - i     a_0 \right)^2,  & i=1,\ldots,M.
  \end{cases}
\end{split}
\end{equation}
Here $a_0\in\R$ denotes the equilibrium distance, and the moduli $k_0$, $k_1$,
and $k_2$ describe the strength of the misfit energy and the elastic
interactions, respectively.  To ensure coercivity, we require $k_0>0$ and
$k_1+2k_2>2|k_2|$.

Next, the continuum energy
\begin{align}
  \E^c_i = \E^{c,e}_i + \E^{c,m}_i
\end{align}
for each atom $i$ is derived from the atomistic energy, which leads to
\begin{equation} \label{EqDefEci}
\begin{split}
  \E^{c,e}_i(\vec y) & =
              \fourth k_{12}(y_i     - y_{i-1} - a_0)^2
            + \fourth k_{12}(y_{i+1} - y_i     - a_0)^2, \\
  \E^{c,m}_i(\vec y) & = \begin{cases}
    \half k_0 \left( y_i - (i-1) a_0 \right)^2,  & i=-M+1,\ldots,0, \\
    \half k_0 \left( y_i - i     a_0 \right)^2,  & i=1,\ldots,M,
  \end{cases}
\end{split}
\end{equation}
where $k_{12}=k_1+4k_2$, see \cite{ArndtLuskin:2007a} for the details.
We note that $\E^{c,e}_i(\vec y)= \E^{a,e}_i(\vec y)$ if
$y_{i+2}-y_{i+1}=y_{i+1}-y_{i}=y_{i}-y_{i-1}=y_{i-1}-y_{i-2},$
but $\E^{c,e}_i(\vec y)\ne \E^{a,e}_i(\vec y)$ more generally.

For each atom $i$, we decide whether this atom is modeled atomistically or as
continuum. An {\em a posteriori} error estimator for this task has been derived
in \cite{ArndtLuskin:2007a}. Let
\begin{align}
\delta_i^a = \begin{cases}
  1 & \text{if atom $i$ is modeled atomistically,} \\
  0 & \text{if atom $i$ is modeled as continuum,}
\end{cases}
\qquad \text{and} \qquad \delta_i^c = 1-\delta_i^a.
\end{align}
We define the {\em atomistic-continuum energy} to be
\begin{equation} \label{EqDefEac}
  \E^{ac}(\vec y)
  :=   \sum_{i=-M+1}^M \delta^a_i \E^a_i(\vec y)
     + \sum_{i=-M+1}^M \delta^c_i \E^c_i(\vec y).
\end{equation}

The next step is to coarsen out unnecessary atoms within the continuum region.
This way, we restrict the system to the remaining atoms, called the
representative atoms, or briefly repatoms. An {\em a posteriori} error estimator
to determine the optimal coarsening will be developed in this paper.  Let
$\ell_j$ for $j=-N+1,\ldots, N$ be the index of the $j$-th repatom, where $2N$
is the number of repatoms.  We require that
\begin{align}
  \ell_{-N+1} < \ell_{-N+2} < \cdots < \ell_j <
  \ell_{j+1} < \cdots < \ell_{N-1} < \ell_N
\end{align}
and
\begin{align}
  \ell_{-N+1} & = -M+1, &
  \ell_{-N+2} & = -M+2, &
  \ell_{N-1} & = M-1, &
  \ell_N & = M.
\end{align}
Then
\begin{align}
  \nu_j := \ell_{j+1}-\ell_j
\end{align}
gives the number of atomistic intervals between the repatoms $j$ and $j+1$. We
denote the vector of repatoms by $\vec{y}^{qc}$.

The {\em quasicontinuum energy} is obtained by implicitly reconstructing the
missing atoms from the repatoms by piecewise linear interpolation, and then
computing the atomistic-continuum energy from the reconstructed vector. The
piecewise linear interpolation can be written as the matrix multiplication
$I^{aq} \vec{y}^{qc}$, where
\begin{align}
  I^{aq}_{\ell_j+k,j  } := \frac{\nu_j-k}{\nu_j}, \qquad
  I^{aq}_{\ell_j+k,j+1} := \frac{      k}{\nu_j}
\end{align}
for $j=-N+1,\ldots,N$ and $k=0,\ldots,\nu_j$, and $I^{aq}_{i,j} = 0$ otherwise.
Hence, the quasicontinuum energy is given by
\begin{align}
  \E^{qc}(\vec y^{qc}) := \E^{ac}(I^{aq} \vec y^{qc}).
\end{align}
Summation formulas for an efficient computation of $\E^{qc}(\vec y^{qc})$
without having to explicitly reconstruct the non-repatoms have been derived in
\cite{ArndtLuskin:2007a}.

We describe boundary conditions at the two leftmost atoms and the two rightmost
atoms.  We take into account two instead of one atom at each end of the chain
because of the NNN interactions.  For the fully atomistic system and the
atomistic-continuum system, the boundary conditions read as
\begin{equation}
\begin{aligned}
  y^{a }_{-M+1} & = y^{bc}_{l1}, & \qquad
  y^{a }_{-M+2} & = y^{bc}_{l2}, & \qquad
  y^{a }_{M-1}  & = y^{bc}_{r2}, & \qquad
  y^{a }_{M}    & = y^{bc}_{r1},  \\
  y^{ac}_{-M+1} & = y^{bc}_{l1}, &
  y^{ac}_{-M+2} & = y^{bc}_{l2}, &
  y^{ac}_{M-1}  & = y^{bc}_{r2}, &
  y^{ac}_{M}    & = y^{bc}_{r1}, \\
\end{aligned}
\end{equation}
for given values $y^{bc}_{l1}$, $y^{bc}_{l2}$, $y^{bc}_{r2}$, and $y^{bc}_{r1}$,
and similarly the boundary conditions for the quasicontinuum system are given by
\begin{equation}
\begin{aligned}
  y^{qc}_{-N+1} & = y^{bc}_{l1}, & \qquad
  y^{qc}_{-N+2} & = y^{bc}_{l2}, & \qquad
  y^{qc}_{N-1}  & = y^{bc}_{r2}, & \qquad
  y^{qc}_{N}    & = y^{bc}_{r1}.
\end{aligned}
\end{equation}

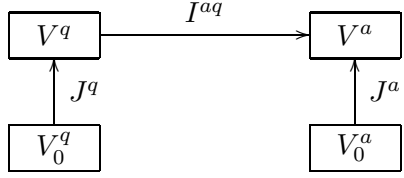
\begin{figure}
\begin{equation*}
\begin{xy}
(  0,15) *\txt{$V^q$}   *=(12,6)\frm{-} ="boxqc";
( 40,15) *\txt{$V^a$}    *=(12,6)\frm{-} ="boxac";
(  0,00) *\txt{$V^q_0$} *=(12,6)\frm{-} ="boxqc0";
( 40,00) *\txt{$V^a_0$}  *=(12,6)\frm{-} ="boxac0";
"boxqc"; "boxac"  **\dir{-} ?>* \dir{>};
? *_!/3mm/{I^{aq}};
"boxqc0"; "boxqc" **\dir{-} ?>* \dir{>};
? *^!/4mm/{J^q};
"boxac0"; "boxac" **\dir{-} ?>* \dir{>};
? *^!/4mm/{J^a};
\end{xy}
\end{equation*}
\caption{Solution spaces and their operators for the fully-atomistic level (a)
  and the quasicontinuum level (q).}
\label{FigAcQc}
\end{figure}

Next, we introduce the solution spaces
\begin{equation}
\begin{aligned}
  V^a   & := \R^{2M},   & \qquad
  V^a_0 & := \R^{2M-4}, & \qquad
  V^q   & := \R^{2N},   & \qquad
  V^q_0 & := \R^{2N-4}
\end{aligned}
\end{equation}
and the boundary operators
\begin{equation}
  J^a: V^a_0 \to V^a \qquad \text{and} \qquad J^q: V^q_0 \to V^q
\end{equation}
by
\begin{equation}
\begin{aligned}
  J^a_{i,j} & := \delta_{ij} &
  \quad \text{for} \quad i & = -M+1,\ldots,M, & j & = -M+3,\ldots,M-2, \\
  J^q_{i,j} & := \delta_{ij} &
  \quad \text{for} \quad i & = -N+1,\ldots,N, & j & = -N+3,\ldots,N-2.
\end{aligned}
\end{equation}
$V^a$ and $V^a_0$ are the suitable spaces for the atomistic system and the
atomistic-continuum system, both with and without boundary values, whereas $V^q$
and $V^q_0$ are the suitable spaces for the quasicontinuum system. $J^a$ and
$J^q$ extend vectors from $V^a_0$ and $V^q_0$, respectively, by zero boundary
values.  Note that the interpolation operator $I^{aq}$ maps $V^q$ to $V^a$. The
spaces and their operators are illustrated in Figure~\ref{FigAcQc}.

To implement the boundary conditions, we define the vectors
\begin{equation}
\begin{aligned} \label{EqDefYbc}
  \vec y^{bcq} & := \begin{bmatrix}
    y^{bc}_{l1} & y^{bc}_{l2} & 0 & \cdots & 0 &
    y^{bc}_{r2} & y^{bc}_{r1}
  \end{bmatrix}^T \in V^q, \\
  \vec y^{bca} & := I^{aq} \vec y^{bcq} \in V^a.
\end{aligned}
\end{equation}

We seek minimizers of the three different potential energies subject to the
boundary conditions described above, that is, vectors
\begin{align}
  \vec {\bar y}^a    \in J^a V^a_0 + \vec y^{bca} \subset V^a,
  \quad
  \vec {\bar y}^{ac} \in J^a V^a_0 + \vec y^{bca} \subset V^a,
  \quad \text{and} \quad
  \vec {\bar y}^{qc} \in J^q V^q_0 + \vec y^{bcq} \subset V^q
\end{align}
which minimize the potential energies $\E^a$, $\E^{ac}$, and $\E^{qc}$,
respectively.  If we decompose
\begin{align}
  \vec {\bar y}^{a}  = J^a \vec y^{a}  + \vec y^{bca}, \qquad
  \vec {\bar y}^{ac} = J^a \vec y^{ac} + \vec y^{bca}, \qquad \text{and} \qquad
  \vec {\bar y}^{qc} = J^q \vec y^{qc} + \vec y^{bcq},
\end{align}
then the fully atomistic solution $\vec y^{a} \in V^a_0$, the
atomistic-continuum solution $\vec y^{ac} \in V^a_0,$ and the quasicontinuum
solution $\vec y^{qc} \in V^q_0$ are characterized by
\begin{equation}
\begin{aligned}
  \vec y^{a}  & = \argmin_{\vec y\in V^a_0} \E^{a} (J^a \vec y+{\vec y}^{bca}), \\
  \vec y^{ac} & = \argmin_{\vec y\in V^a_0} \E^{ac}(J^a \vec y+{\vec y}^{bca}), \\
  \vec y^{qc} & = \argmin_{\vec y\in V^q_0} \E^{qc}(J^q \vec y+{\vec y}^{bcq}).
\end{aligned}
\end{equation}

Next, we write the energies in matrix notation as
\begin{equation}
\begin{aligned}
\E^{a }(\vec y) & =   \half (      \vec y-\vec a^a)^T D^T E^{a } D
                            (      \vec y-\vec a^a)
                    + \half (      \vec y-\vec b^a)^T     K
                            (      \vec y-\vec b^a),               \\
\E^{ac}(\vec y) & =   \half (      \vec y-\vec a^a)^T D^T E^{ac} D
                            (      \vec y-\vec a^a)
                    + \half (      \vec y-\vec b^a)^T     K
                            (      \vec y-\vec b^a),               \\
\E^{qc}(\vec y) & =   \half (I^{aq}\vec y-\vec a^a)^T D^T E^{ac} D
                            (I^{aq}\vec y-\vec a^a)
                    + \half (I^{aq}\vec y-\vec b^a)^T     K
                            (I^{aq}\vec y-\vec b^a).
\end{aligned}
\end{equation}
We refer to \cite{ArndtLuskin:2007a} for the precise and lengthy definitions of
the respective matrices. In matrix notation, the vectors $\vec y^a$, $\vec
y^{ac}$, and $\vec y^{qc}$ are the solutions of the linear systems
\begin{equation} \label{EqPrimal}
\begin{aligned}
  M^{a } \vec y^{a } & = \vec f^{a }, \\
  M^{ac} \vec y^{ac} & = \vec f^{ac}, \\
  M^{qc} \vec y^{qc} & = \vec f^{qc},
\end{aligned}
\end{equation}
where the matrices $M^a$, $M^{ac}$, and $M^{qc}$ are given by
\begin{equation} \label{EqDefM}
\begin{split}
M^{a }      & :=   J^{a^T}          ( D^T E^{a } D + K)        J^a, \\
M^{ac}      & :=   J^{a^T}          ( D^T E^{ac} D + K)        J^a, \\
M^{qc}      & :=   J^{q^T} I^{aq^T} ( D^T E^{ac} D + K) I^{aq} J^q, \\
\end{split}
\end{equation}
and where the right-hand sides $\vec f^a$, $\vec f^{ac}$, and $\vec f^{qc}$
are defined as
\begin{equation} \label{EqDefRHS}
\begin{split}
\vec f^{a } & := - J^{a^T}          D^T E^{a}  D (       \vec y^{bca} - \vec a^a)
                 - J^{a^T}              K        (       \vec y^{bca} - \vec b^a), \\
\vec f^{ac} & := - J^{a^T}          D^T E^{ac} D (       \vec y^{bca} - \vec a^a)
                 - J^{a^T}              K        (       \vec y^{bca} - \vec b^a), \\
\vec f^{qc} & := - J^{q^T} I^{aq^T} D^T E^{ac} D (I^{aq} \vec y^{bcq} - \vec a^a)
                 - J^{q^T} I^{aq^T}     K        (I^{aq} \vec y^{bcq} - \vec b^a).
\end{split}
\end{equation}


\section{Goal-Oriented Error Estimation}

We estimate the error $\vec y^{a} - J^{a^T} I^{aq} J^q \vec
y^{qc}$ in terms of a user-definable goal function
\begin{equation}
  Q: V^a_0 \to \R,
\end{equation}
that is, we aim at estimating
\begin{equation}
  |Q(\vec y^{ac}) - Q(J^{a^T} I^{aq} J^q \vec y^{qc})|.
\end{equation}
We assume that $Q$ is linear. Hence there exists some vector $\vec q\in
V^a_0$ such that
\begin{align} \label{EqErrorAQc}
  Q(\vec y) = \vec q^T \vec y \qquad \forall \vec y \in V^a_0.
\end{align}

We decompose the error \eqref{EqErrorAQc} as follows:
\begin{equation} \label{EqErrorPrimal}
\begin{aligned}
  |Q(\vec y^a) - Q( J^{a^T} I^{aq} J^q \vec y^{qc})|
  & =   |Q(\vec y^a - \vec y^{ac})
      +  Q(\vec y^{ac} - J^{a^T} I^{aq} J^q \vec y^{qc})| \\
  & \le |Q(\vec y^a - \vec y^{ac})| 
      + |Q(\vec y^{ac} - J^{a^T} I^{aq} J^q \vec y^{qc})| \\
  & =   |Q(\vec e^{a-ac})| + |Q(\vec e^{ac-qc})|
\end{aligned}
\end{equation}
where
\begin{equation}
\begin{aligned}
  \vec e^{ a-ac} & := \vec y^a - \vec y^{ac}, \\
  \vec e^{ac-qc} & := \vec y^{ac} - J^{a^T} I^{aq} J^q \vec y^{qc}.
\end{aligned}
\end{equation}
The first term $|Q(\vec e^{a-ac})|$ constitutes the modeling error and has been
treated in \cite{ArndtLuskin:2007a}. The second term $|Q(\vec e^{ac-qc})|$
describes the error due to mesh coarsening and will be treated in the following.

To facilitate the error analysis, we define the dual problems
\begin{equation}
\begin{aligned} \label{EqDual}
  M^{ac} \vec g^{ac} & =                      \vec q, \\
  M^{qc} \vec g^{qc} & = J^{q^T} I^{aq^T} J^a \vec q.
\end{aligned}
\end{equation}
We then have the basic dual identity for the goal-oriented error
\begin{equation} \label{EqBasicDual}
\begin{split} 
  Q(\vec e^{ac-qc}) 
  & = \vec q^T (\vec y^{ac} - J^{a^T} I^{aq} J^q \vec y^{qc}) \\
  & = \vec g^{ac^T} M^{ac} (\vec y^{ac} - J^{a^T} I^{aq} J^q \vec y^{qc}) \\
  & = \vec g^{ac^T} \residual^{ac}(J^{a^T} I^{aq} J^q \vec y^{qc})
\end{split}
\end{equation}
with the primal residual given by
\begin{align}
  \residual^{ac}(\vec y) := \vec f^{ac} - M^{ac} \vec y.
\end{align}

However, this quantity is too expensive to compute. We cannot solve for the dual
solution $\vec g^{ac}$ on the atomistic scale, since this has the same
computational complexity as solving for the uncoarsened solution $\vec y^{ac}$.
To overcome this obstacle, we replace $\vec g^{ac}$ by a dual solution from a
coarser space.

A first idea would be to use $J^{a^T} I^{aq} J^q \vec g^{qc}$ instead of $\vec
g^{ac}$. However, this turns out to be useless due to Galerkin orthogonality:

\begin{lemma}[Garlerkin Orthogonality] \label{LemmaGalerkinOrth}
  We have
  \begin{equation}
    J^{q^T} I^{aq^T} J^a M^{ac} ( \vec y^{ac} - J^{a^T} I^{aq} J^q \vec y^{qc} ) = 0,
  \end{equation}
  or equivalently
  \begin{equation}
    J^{q^T} I^{aq^T} J^a \residual^{ac} ( J^{a^T} I^{aq} J^q \vec y^{qc} ) = 0.
  \end{equation}
\end{lemma}

\begin{proof}
  Multiplying the equation for $\vec y^{ac}$ from \eqref{EqPrimal} by $J^{q^T}
  I^{aq^T} J^a$ from the left gives
  \begin{equation}  \label{EqGOProof1}
    J^{q^T} I^{aq^T} J^a M^{ac} \vec y^{ac}
    = J^{q^T} I^{aq^T} J^a \vec f^{ac}.
  \end{equation}
  It is easy to see that
  \begin{equation} \label{EqGOProof3}
    J^a J^{a^T} I^{aq} J^q  = I^{aq} J^q.
  \end{equation}
  Applying this to the equation for $\vec y^{qc}$ from \eqref{EqPrimal} leads to
  \begin{equation} \label{EqGOProof2}
    J^{q^T} I^{aq^T} J^a M^{ac} J^{a^T} I^{aq} J^q \vec y^{qc} = \vec f^{qc}.
  \end{equation}
  Subtracting \eqref{EqGOProof2} from \eqref{EqGOProof1}, substituting the
  definitions \eqref{EqDefRHS} of $\vec f^{ac}$ and $\vec f^{qc}$ and using
  \eqref{EqGOProof3} and \eqref{EqDefYbc} gives
  \begin{equation}
  \begin{split}
    J^{q^T} I^{aq^T} & J^a M^{ac} ( \vec y^{ac} - J^{a^T} I^{aq} J^q \vec y^{qc} ) \\
    & = J^{q^T} I^{aq^T} J^a \vec f^{ac} - \vec f^{qc} \\
    & = J^{q^T} I^{aq^T} J^a J^{a^T}
        \left[ - D^T E^{ac} D (\vec y^{bc} - \vec a^a)
               -     K        (\vec y^{bc} - \vec b^a) \right] \\
    & \qquad - J^{q^T} I^{aq^T}
        \left[ - D^T E^{ac} D ( I^{aq} \vec y^{bcq} - \vec a^a)
               -     K        ( I^{aq} \vec y^{bcq} - \vec b^a) \right] \\
    & = 0,
  \end{split}
  \end{equation}
  which completes the proof.
\end{proof}

Hence, replacing $\vec g^{ac}$ in \eqref{EqBasicDual} by $J^{a^T} I^{aq} J^q \vec
g^{qc}$ always gives a zero estimate for the goal-oriented error. We need to use the dual solution
from some space which is finer than $V^q_0$ to get a non-zero estimate for
the goal-oriented error,
but which is coarser than $V^a_0$ to make it computable.

\begin{figure}
\includegraphics[width=0.8\textwidth]{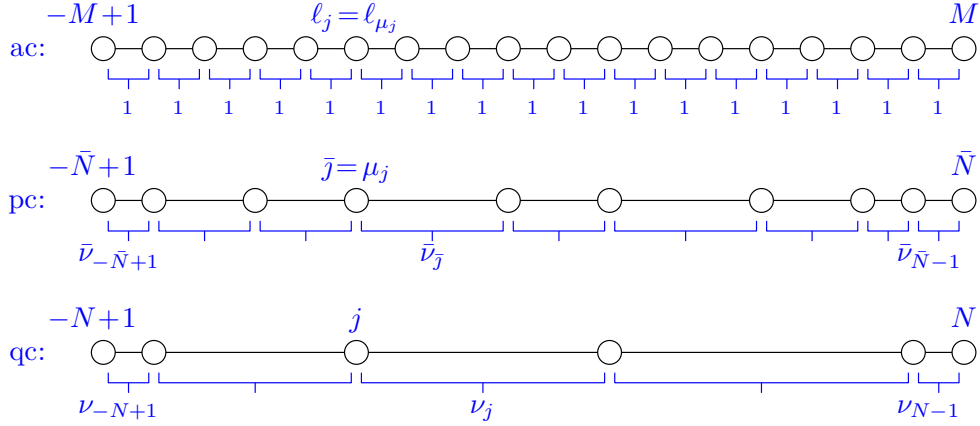}
\caption{Levels of refinement: atomistic-continuum (ac), partial continuum (pc),
  and quasicontinuum (qc) with indices (above the respective chain) and
  intervals (below the respective chain). }
\label{FigPcQcRefinement}
\end{figure}

To this end, we introduce an additional level of refinement and denote it as the
{\em partial continuum} (pc) level.  Similarly to the repatoms on the qc-level,
we chose a set of $2\bar N$ pc-level repatoms which is a subset of all atoms and
a superset of the qc-repatoms.  This means we choose indices $\bar\ell_{\bar\jmath}$
for $\bar\jmath=-\bar N+1, \ldots, \bar N$ such that
\begin{align}
  \bar\ell_{-\bar N+1} < \bar\ell_{-\bar N+2} < \cdots < \bar\ell_{\bar\jmath} <
  \bar\ell_{\bar\jmath+1} < \cdots < \bar\ell_{\bar N-1} < \bar\ell_{\bar N}
\end{align}
and
\begin{align}
  \bar\ell_{-\bar N+1} & = -M+1, &
  \bar\ell_{-\bar N+2} & = -M+2, &
  \bar\ell_{ \bar N-1} & = M-1, &
  \bar\ell_{ \bar N  } & = M.
\end{align}
To ensure that the pc-level is actually a refinement of the qc-level, we require
every qc-level repatom to be a pc-level repatom. Hence there exist indices
$\mu_j$ such that
\begin{equation}
 \ell_j = \bar\ell_{\mu_j}, \qquad j=-N+1,\ldots, N.
\end{equation}
Similar to the definition of $\nu_j$, we denote the number of atomistic
intervals between two pc-level repatoms by
\begin{equation}
  \bar\nu_{\bar\jmath} := \bar\ell_{\bar\jmath+1} - \bar\ell_{\bar\jmath}, \qquad
  \bar\jmath = -\bar N+1,\ldots,\bar N-1.
\end{equation}
See Figure~\ref{FigPcQcRefinement} for an illustration of the three levels of
refinement and the corresponding variables.

\begin{figure}
\begin{align}
\begin{xy}
(  0,15) *\txt{$V^q$}   *=(12,6)\frm{-} ="boxqc";
( 40,15) *\txt{$V^p$}   *=(12,6)\frm{-} ="boxpc";
( 80,15) *\txt{$V^a$}    *=(12,6)\frm{-} ="boxac";
(  0,00) *\txt{$V^q_0$} *=(12,6)\frm{-} ="boxqc0";
( 40,00) *\txt{$V^p_0$} *=(12,6)\frm{-} ="boxpc0";
( 80,00) *\txt{$V^a_0$}  *=(12,6)\frm{-} ="boxac0";
"boxqc"!<0em,0.2em>; "boxpc"!<0em,0.2em> **\dir{-} ?>* \dir{>};
? *_!/0.5em/{I^{pq}};
"boxpc"!<0em,-0.2em>; "boxqc"!<0em,-0.2em> **\dir{-} ?>* \dir{>};
? *_!/0.5em/{R^{qp}};
"boxpc"; "boxac" **\dir{-} ?>* \dir{>};
? *_!/0.5em/{I^{ap}};
"boxqc"; "boxac" **\crv{(40,30)} ?>* \dir{>};
? *_!/0.5em/{I^{aq}=I^{ap}I^{pq}};
"boxqc0"; "boxqc" **\dir{-} ?>* \dir{>};
? *^!/0.7em/{J^q};
"boxpc0"; "boxpc" **\dir{-} ?>* \dir{>};
? *^!/0.7em/{J^p};
"boxac0"; "boxac" **\dir{-} ?>* \dir{>};
? *^!/0.7em/{J^a};
\end{xy}
\end{align}
\caption{Solution spaces and their operators for the ac-level, pc-level, and
  qc-level.}
\label{FigAcPcQc}
\end{figure}
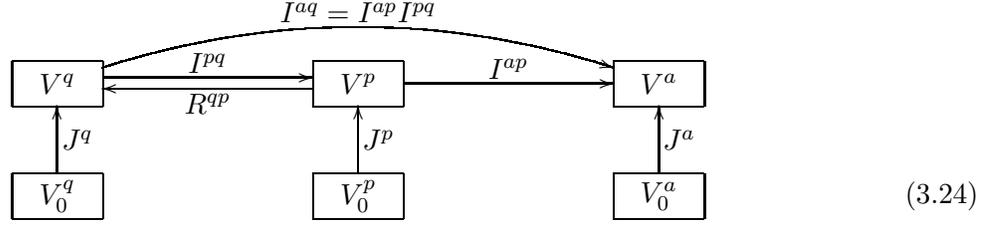

We define the pc-level solution spaces
\begin{equation}
  V^p   := \R^{2\bar N} \qquad \text{and} \qquad
  V^p_0 := \R^{2\bar N-4},
\end{equation}
the interpolation operators
\begin{equation}
\begin{aligned}
  \left. \begin{aligned}
    I^{ap}_{\bar\ell_{\bar\jmath}+\bar k,\bar\jmath  } &
      := \makebox[7em][l]{$\displaystyle \frac{\bar\nu_{\bar\jmath}-\bar k}{\bar\nu_{\bar\jmath}}$} \\
    I^{ap}_{\bar\ell_{\bar\jmath}+\bar k,\bar\jmath+1} &
      := \makebox[7em][l]{$\displaystyle \frac{                     \bar k}{\bar\nu_{\bar\jmath}}$}
  \end{aligned} \quad \right\} &
  \quad \text{for } \bar\jmath = -\bar N+1,\ldots,\bar N-1,
              \quad \bar k = 0,\ldots,\bar\nu_{\bar\jmath}, \\
  I^{ap}_{i,\bar\jmath}
      := \makebox[8.9em][l]{$0$} & \quad \text{otherwise},\\
  \left. \begin{aligned}
    I^{pq}_{\mu_j+k,j  } &
      := \makebox[7em][l]{$\displaystyle \frac{\nu_j-\bar\ell_{\mu_j+k}+\bar\ell_{\mu_j}}{\nu_j}$} \\
    I^{pq}_{\mu_j+k,j+1} &
      := \makebox[7em][l]{$\displaystyle \frac{\bar\ell_{\mu_j+k}-\bar\ell_{\mu_j}}{\nu_j}$}
  \end{aligned} \quad \right\} &
  \quad \text{for } j = -N+1,\ldots,N-1,
              \quad k = 0,\ldots,\mu_{j+1}-\mu_j, \\
  I^{pq}_{\bar\jmath,j}
      := \makebox[8.9em][l]{$0$} & \quad \text{otherwise},
\end{aligned}
\end{equation}
the restriction operator
\begin{align}
  R^{qp}_{j,\bar\jmath} & := \delta_{\mu_j,\bar\jmath}
  \qquad \text{for} \quad j = -N+1,\ldots,N,
                    \quad \bar\jmath = -\bar{N}+1,\ldots,\bar{N},
\end{align}
and the boundary operator
\begin{align}
  J^p_{i,j} & := \delta_{ij}
  \qquad \text{for} \quad i = -\bar{N}+1,\ldots,\bar{N},
                    \quad j = -\bar{N}+3,\ldots,\bar{N}-2.
\end{align}
Note that the interpolation operators are defined in such a way that
$I^{aq}$ factors as
\begin{align}
  I^{aq} = I^{ap} I^{pq}.
\end{align}
Altogether, all solution spaces and the corresponding operators
are depicted in Figure~\ref{FigAcPcQc}.

We solve the dual equation on the newly introduced partial continuum level:
\begin{align}
  M^{pc} \vec g^{pc} & = J^{p^T} I^{ap^T} J^a \vec q,   \label{EqDualPc}
\end{align}
where
\begin{align}
  M^{pc} & := J^{p^T} I^{ap^T} ( D^T E^{ac} D + K) I^{ap} J^p.
\end{align}
We use $\vec g^{pc}$ to get an approximation of the goal-oriented error:
\begin{align} \label{EqErrorApprox}
  Q(\vec y^{ac}) - Q(J^{a^T} I^{aq} J^q \vec y^{qc})
  \approx \vec g^{pc^T} J^{p^T} I^{ap^T} J^a \residual^{ac}(J^{a^T} I^{aq} J^q \vec y^{qc}).
\end{align}

Due to Galerkin Orthogonality (Lemma~\ref{LemmaGalerkinOrth}), we can subtract
any vector in $V^q_0$ from the right hand side expression without changing its
value. Thus, we subtract from $\vec g^{pc}$ the linear interpolation of the
nodal values of $\vec g^{pc}$ at the qc-repatoms, $J^{p^T} I^{pq} R^{qp} J^p
\vec g^{pc}$, to get a vector in $V^p_0$ that is zero at the qc-repatoms.  This
leads us to the error estimator
\begin{equation} \label{EqEE}
\begin{split}
  \eta:=&\:\big[ \vec g^{pc} - J^{p^T} I^{pq} R^{qp} J^p \vec g^{pc} \big]^T
           J^{p^T} I^{aq^T} J^a \residual^{ac}(J^{a^T} I^{aq} J^q \vec y^{qc}) \\
       =&\:\big[ \vec g^{pc} - J^{p^T} I^{pq} R^{qp} J^p \vec g^{pc} \big]^T
           \big[ J^{p^T} I^{aq^T} J^a \vec f^{ac} - M^{pc} \vec y^{pc} \big]
\end{split}
\end{equation}
as an approximation to \eqref{EqBasicDual}.

In Section~\ref{SecResults}, we will construct a numerical method based on this
error estimator. We will then compute numerically how much the approximation
\eqref{EqEE} deviates from the precise error \eqref{EqBasicDual}, and we will
determine how fine the partial refinement needs to be.


\section{Numerical Method}

Based on the error estimator $\eta$, we need to decide what
intervals at the qc-level shall be refined.  To this end, we split
$\eta$ into individual contributions from each pc-level repatom:
\begin{align} \label{EqEtaPcSplit}
  \eta = \sum_{\bar\jmath=-\bar N+3}^{\bar N-2} \eta^{pc}_{\bar\jmath}
\end{align}
where
\begin{align} \label{EqEtaSplit1}
  \eta^{pc}_{\bar\jmath} :=
    \big[ \vec g^{pc} - J^{p^T} I^{pq} R^{qp} J^p \vec g^{pc} \big]_{\bar\jmath}
    \big[ J^{p^T} I^{aq^T} J^a \residual^{ac}(J^{a^T} I^{aq} J^q \vec y^{qc}) \big]_{\bar\jmath}.
\end{align}
Next, we map the values $\eta^{pc}_{\bar\jmath}$ to the qc-level intervals.
Note that we already have $\eta^{pc}_{\bar\jmath}=0$ whenever $\bar\jmath$ is a
qc-repatom because we subtracted the qc-level interpolant from $\vec g^{pc}$ in
\eqref{EqEE}. We accumulate the remaining values $\eta^{pc}_{\bar\jmath}$ to the
interval between the qc-repatoms $j$ and $j+1$ if the pc-repatom ${\bar\jmath}$
lies between the qc-repatoms $j$ and $j+1$, that is, $\mu_j < {\bar\jmath} <
\mu_{j+1},$ to obtain
\begin{equation} \label{EqEtaSplit2}
  \eta^{qc}_j := \Bigg| \sum_{\bar\jmath=\mu_j+1}^{\mu_{j+1}-1}
  \eta^{pc}_{\bar\jmath} \Bigg|.
\end{equation}

Moreover, we have to decide how to compute the partial refinement. Here we
subdivide each qc-interval $(\ell_j, \ell_{j+1})$ into a fixed number $\Lambda$
of subintervals of the same size, although other strategies are possible, like
letting $\Lambda$ depend on the respective interval size $\nu_j = \ell_{j+1} -
\ell_j$, or dividing each interval into subintervals of different sizes. But
here we stick to a fixed number $\Lambda\in\N$ and an equidistant
subdivision. We employ the following algorithm for partial refinement:

\noindent
\hspace*{20mm} $\bar\jmath \leftarrow -\bar{N}+1$\\
\hspace*{20mm} for $j=-N+1,\ldots,N-1$\\
\hspace*{30mm}     $\omega \leftarrow \max(1, \nu_j/\Lambda$)\\
\hspace*{30mm}     $\sigma_1 \leftarrow 0$\\
\hspace*{30mm}     $\sigma_2 \leftarrow 0$\\
\hspace*{30mm}     while $(\sigma_2 < \nu_j)$\\
\hspace*{40mm}         $\sigma_1 \leftarrow \min(\sigma_1+\omega, \nu_j)$\\
\hspace*{40mm}         $\bar{\nu}_{\bar\jmath} \leftarrow \lfloor\sigma_1-\sigma_2+\half \rfloor$\\
\hspace*{40mm}         $\sigma_2 \leftarrow \sigma_2+\bar\nu_{\bar\jmath}$\\
\hspace*{40mm}         $\bar\jmath \leftarrow \bar\jmath+1$\\
\hspace*{30mm}     end\\
\hspace*{20mm} end

Note that this algorithm performs the necessary rounding if $\nu_j$ is not
divisible by $\Lambda$. If $\nu_j \le \Lambda$, then the interval gets fully
refined up to the atomistic level.

We then use this partial refinement algorithm and the global and local error
estimators above to construct the following algorithm for adaptive mesh
coarsening:

\begin{enumerate}
\item[(1)] Start with the model fully coarsened in the continuum region.
\item[(2)] Solve the primal problem \eqref{EqPrimal} for $\vec y^{qc}.$
\item[(3)] Determine the partial refinement pc according to the algorithm above.
\item[(4)] Solve the dual problem \eqref{EqDualPc} for $\vec g^{pc}$.
\item[(5)] Compute the global error estimator $\eta$ from \eqref{EqEE}.
\item[(6)] If $\eta \le \tau_{gl}$, then stop.
\item[(7)] Compute the local error estimators $\eta^{qc}_j$ from
  \eqref{EqEtaSplit1} and \eqref{EqEtaSplit2}.
\item[(8)] Refine all intervals $(j,j+1)$ with
  \begin{align} \label{EqAdaptCrit}
    \eta^{qc}_j \ge \frac{1}{\tau_{fac}} \max_k \eta^{qc}_k
  \end{align}
  into two subintervals.
\item[(8)] Go to (2).
\end{enumerate}

Here $\tau_{gl}>0$ denotes the global error tolerance.  For the adaption
criterion \eqref{EqAdaptCrit}, numerical evidence has shown that $\tau_{fac}=10$
is a reasonable choice.

A possible improvement of step (8) would be to use the absolute value of the
error estimator $\eta^{qc}_j$ to decide into how many subintervals each interval
$(j,j+1)$ shall be refined, instead of just refining into two
subintervals. However, the magnitudes of the values $\eta^{qc}_j$ differ
considerably only in the first iterations.  So this would only eliminate some of
the early iterations which are still cheap due to the small number of unknowns
involved there, whereas the computationally expensive later iterations will not
be affected. Hence we expect the overall speedup to be small.


\section{Numerical Results} \label{SecResults}

Now we present and discuss our numerical results.  We choose
boundary conditions
\begin{equation}
  y^{bc}_{l1} = -M,   \qquad
  y^{bc}_{l2} = -M+1, \qquad
  y^{bc}_{r2} = M-1,  \qquad
  y^{bc}_{r1} = M.
\end{equation}
The elastic moduli are given by $k_0=0.1$, $k_1=2$, and $k_2=1$, and the lattice
constant is $a_0=1$.

\begin{figure}
\includegraphics[width=0.7\textwidth]{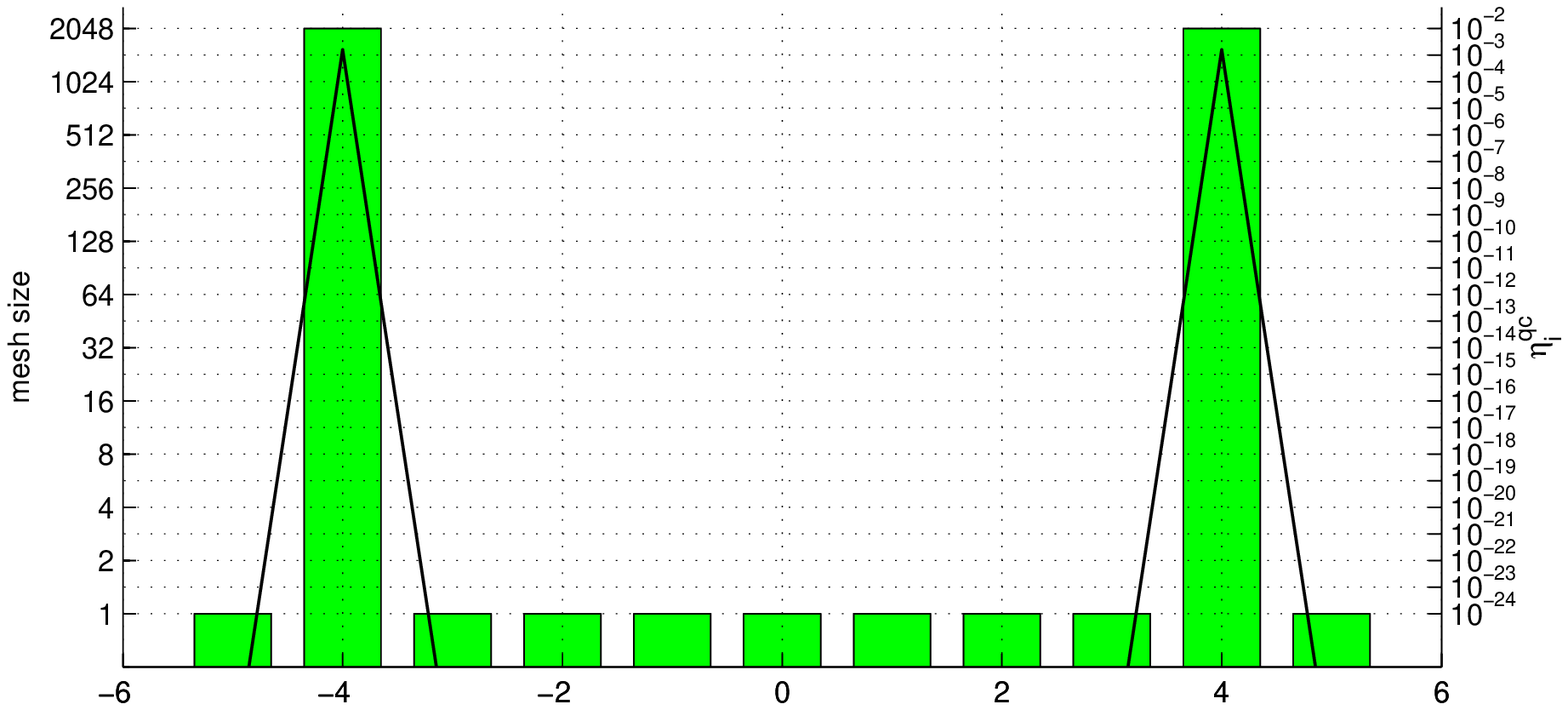}
\includegraphics[width=0.7\textwidth]{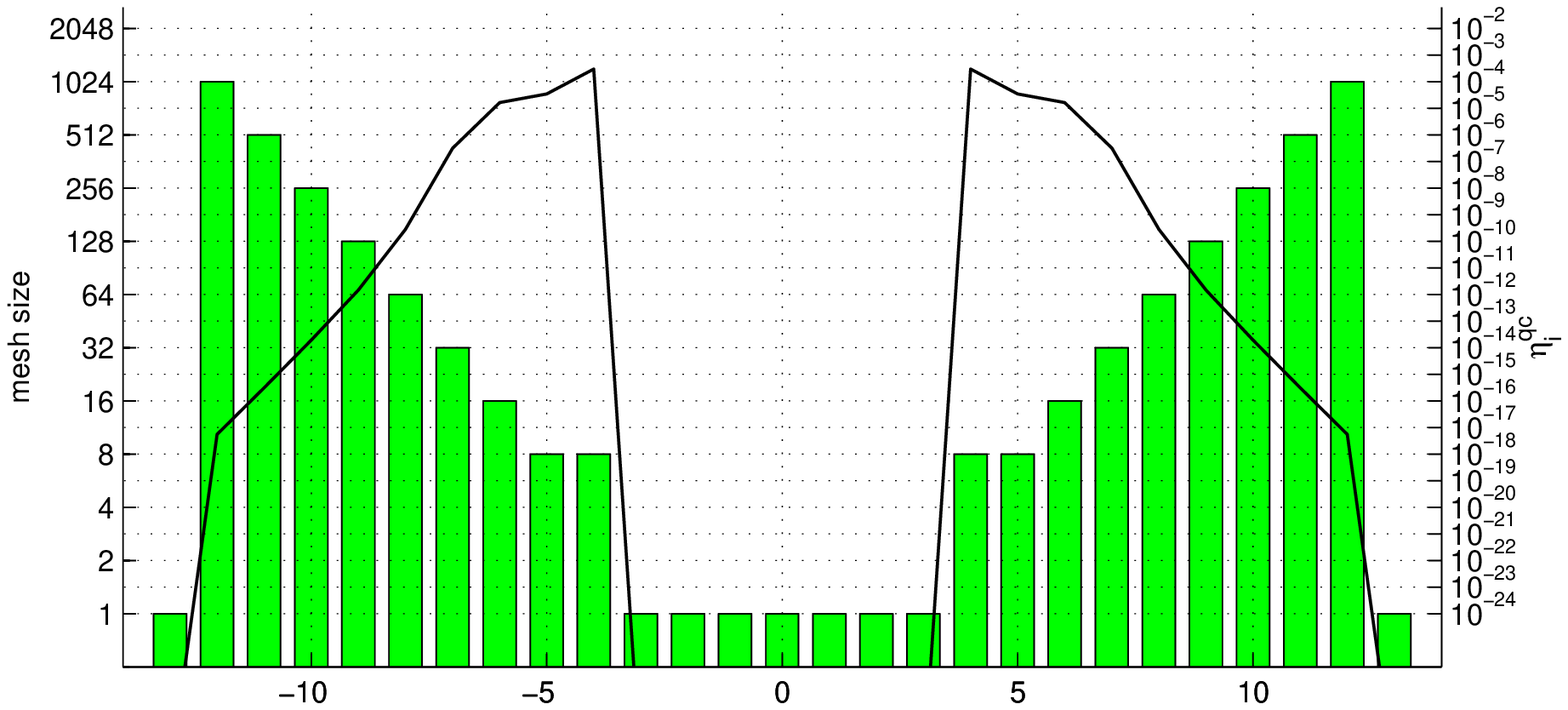}
\includegraphics[width=0.7\textwidth]{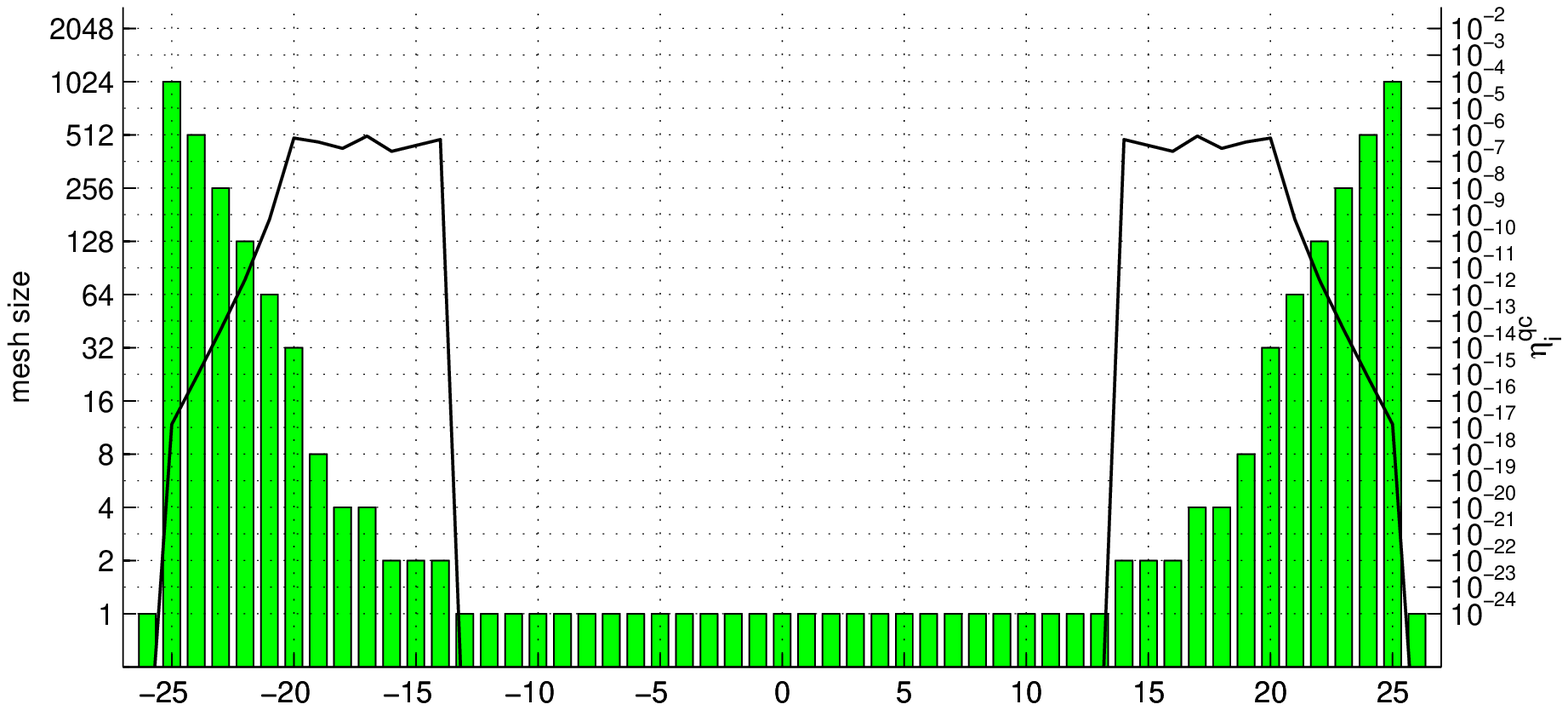}
\caption{Mesh size $\nu_j$ (green bar graph) and $\eta^{qc}_j$ (black line graph)
  for error tolerances $\tau_{gl}=10^{-1}$ (top), $\tau_{gl}=10^{-3}$ (middle)
  and $\tau_{gl}=10^{-5}$ (bottom).}
\label{FigMeshEta}
\end{figure}

We consider a chain of 4106 atoms, that is $M=2053$. There is an atomistic
region of $4$ atoms from $-1$ to $2$ around the dislocation at the center of the
chain. The remaining part is modeled as continuum. We set atoms $-3,\, -2,\,3,
\,4$ to be continuum repatoms so that $\E^a_{-1}(I^{aq} \vec y^{qc})$ and
$\E^a_{2}(I^{aq} \vec y^{qc})$ can be evaluated without interpolation, and we
set and atoms $-M+1,\, -M+2, \,M-1, \,M$ to be continuum repatoms so that the
boundary conditions can be set.  Initially, the mesh in the continuum region is
maximally coarsened, so that there are no more repatoms in addition to the ones
already mentioned.  This gives $N=6$ and 12 qc-level degrees of freedom.  There
are two large elements of size $\nu_{-4} = \nu_4 = 2048$, one on each side of
the dislocation at the center of the chain, whereas the remaining elements
necessarily have sizes $\nu_j = 1$ for $j=-5, -3, -2, -1, 0, 1, 2, 3, 5$. The
mesh is shown in the upper graph of Figure~\ref{FigMeshEta}.

The quantity of interest here is the size $y_1-y_0$ of the dislocation at the
center of the chain, that is, the distance between the two atoms 0 and 1 to the
left and right of the dislocation. The corresponding vector $\vec q\in V^{a}_0$
reads as
\begin{align}
  \vec q = [0,\ldots,0,-1,1,0,\ldots,0]^T.
\end{align}

\begin{table}
\begin{tabular}{|c|c|c|c|c|c|c|}
\hline
iteration & \#dof & $\min \nu_j$ & $\max \nu_j$ &
$\eta$ & $\sum_{j=-N+1}^{N-1} \eta_j^{qc}$ & $|Q(\vec e^{ac-qc})|$ \\
\hline
 1 & 12 & 2048 & 2048 & 3.143618e-03 & 3.143618e-03 & 6.777614e-02 \\
 2 & 14 & 1024 & 1024 & 5.208032e-03 & 5.443530e-03 & 6.463252e-02 \\
 3 & 16 &  512 & 1024 & 8.771892e-03 & 9.133002e-03 & 5.946329e-02 \\
 4 & 18 &  256 & 1024 & 1.293987e-02 & 1.343519e-02 & 5.074706e-02 \\
 5 & 20 &  128 & 1024 & 1.520764e-02 & 1.565599e-02 & 3.787477e-02 \\
 6 & 22 &   64 & 1024 & 1.267077e-02 & 1.279361e-02 & 2.271288e-02 \\
 7 & 24 &   32 & 1024 & 6.760509e-03 & 6.767672e-03 & 1.004707e-02 \\
 8 & 26 &   16 & 1024 & 2.395699e-03 & 2.395699e-03 & 3.286644e-03 \\
 9 & 28 &    8 & 1024 & 6.933383e-04 & 6.933394e-04 & 9.216477e-04 \\
10 & 32 &    4 & 1024 & 2.061976e-04 & 2.061988e-04 & 2.638938e-04 \\
11 & 40 &    2 & 1024 & 5.841551e-05 & 5.841804e-05 & 6.391755e-05 \\
12 & 54 &    1 & 1024 & 7.567732e-06 & 7.570401e-06 & 9.376934e-06 \\
\hline
\end{tabular}
\caption{Advance of the algorithm until the error tolerance $\tau_{gl}=10^{-5}$ is
  achieved.}
\label{TabAlgorithm}
\end{table}

Table~\ref{TabAlgorithm} shows how the algorithm proceeds for an error tolerance
of $\tau_{gl} = 10^{-5}$. One can easily see how the elements get refined and
the error drops.  Also shown is the precise error $|Q(\vec e^{ac-qc})|$ which is
available for this small model problem.  Note that the column $\min \nu_j$
refers only to the elements in the continuum region which actually can be
coarsened, excluding the padding around the atomistic region and the boundary
layer.

The bar graph in Figure~\ref{FigMeshEta} shows the size of the individual
elements in the meshes before the first iteration and after the error tolerances
$\tau_{gl}=10^{-3}$ and $\tau_{gl}=10^{-5}$ are achieved. One can clearly see
how the elements get successively refined towards the center.

The black lines in Figure~\ref{FigMeshEta} depict the decomposed
error estimator $\eta^{qc}_j$. As a result of the Galerkin
orthogonality, it vanishes whereever the mesh is fully refined up
to the atomistic level. We can as well read off the graphs how the
algorithm tends to distribute the error uniformly over the whole
mesh as the error tolerance is decreased. For $\tau_{gl}=10^{-5}$,
we already achieved a close to flat region from elements
$-20\ldots -13$ and $13\ldots20$.

\begin{table}
\begin{tabular}{|c|c|c|c|c|c|c|}
\hline
       & $|Q(\vec e^{ac-qc})|$
                      & $\Lambda$ & $\eta$   & $\sum \eta^{qc}_j$
                                             & $\eta/|Q(\vec e^{ac-qc})|$
                                             & $\sum \eta^{qc}_j/|Q(\vec e^{ac-qc})|$ \\
\hline
mesh 1 & 6.777614e-02 & 2         & 3.143618e-03 & 3.143618e-03 & 0.046382 & 0.046382 \\
       &              & 4         & 8.351650e-03 & 8.351650e-03 & 0.123224 & 0.123224 \\
       &              & 8         & 1.708501e-02 & 1.708501e-02 & 0.252080 & 0.252080 \\
       &              & $\infty$  & 6.777614e-02 & 6.777614e-02 & 1.000000 & 1.000000 \\
\hline
mesh 2 & 9.216477e-04 & 2         & 6.933383e-04 & 6.933394e-04 & 0.752281 & 0.752282 \\
       &              & 4         & 8.749195e-04 & 8.749264e-04 & 0.949299 & 0.949307 \\
       &              & 8         & 9.208978e-04 & 9.209073e-04 & 0.999186 & 0.999197 \\
       &              & $\infty$  & 9.216477e-04 & 9.216582e-04 & 1.000000 & 1.000011 \\
\hline
mesh 3 & 9.376934e-06 & 2         & 7.567732e-06 & 7.570401e-06 & 0.807058 & 0.807343 \\
       &              & 4         & 9.070422e-06 & 9.085687e-06 & 0.967312 & 0.968940 \\
       &              & 8         & 9.320553e-06 & 9.341074e-06 & 0.993987 & 0.996176 \\
       &              & $\infty$  & 9.376934e-06 & 9.399414e-06 & 1.000000 & 1.002397 \\
\hline
\end{tabular}
\caption{Efficiency of the error estimator, $\eta$.}
\label{TabEfficiency}
\end{table}

Table~\ref{TabEfficiency} shows the efficiency of the error estimator for
different meshes and different values of $\Lambda$.  Meshes 1, 2, and 3 refer to
the meshes displayed in Figure~\ref{FigMeshEta}.  For comparison, we also
include the value $\Lambda=\infty$, which indicates that the pc-level mesh for
the dual solution is fully refined to the atomistic level.  One can read off
from this table that the error indicator gets closer to the precise error if the
mesh gets finer. For the coarse mesh 1, the actual error is considerably
underestimated.  However, this has only little impact on the final mesh since
all elements have to be refined anyhow.  For the finer meshes 2 and 3, $\eta$
gets closer to the precise error $|Q(\vec e^{ac-qc})|$.

Also, we can see from Table~\ref{TabEfficiency} how the choice of $\Lambda$
affects the error estimator. As expected, $\eta$ gets closer to the precise
error the larger $\Lambda$ gets. For computational efficiency, we are interested
in keeping $\Lambda$ small, though. We can read off from the table for the fine
mesh 3 that already the smallest possible value $\Lambda=2$ gives a good
estimate of the error, with a deviation of less than 20\%. For $\Lambda=4$, we
already get a very precise estimate.

The absolute accuracy of the error estimator is important, but even more
important is how well $\eta$ controls the mesh refinement, that is how efficient
the resulting mesh is in terms of reaching a prescribed error tolerance with a
minimal number of degrees of freedom.  We now investigate the influence of the
parameter $\Lambda$ on the mesh quality.

\begin{table}
\begin{tabular}{|r||r|c|c||r|c|c|c|c|}
\hline
   & \multicolumn{3}{c||}{$\Lambda=2$} &
     \multicolumn{2}{c|}{} &
     $\Lambda=4$ & $\Lambda=8$ & $\Lambda=\infty$ \\
it & \#dof & $|Q(\vec e^{ac-qc})|$ & $\eta$
   & \#dof & $|Q(\vec e^{ac-qc})|$ & $\eta$ & $\eta$ & $\eta$    \\
\hline
 1 &  12 & 6.778e-02 & 3.144e-03 &  12 & 6.778e-02 & 8.352e-03 & 1.709e-02 & 6.778e-02 \\
 2 &  14 & 6.463e-02 & 5.208e-03 &  14 & 6.463e-02 & 1.394e-02 & 2.683e-02 & 6.463e-02 \\
 3 &  16 & 5.946e-02 & 8.772e-03 &  16 & 5.946e-02 & 2.166e-02 & 3.680e-02 & 5.946e-02 \\
 4 &  18 & 5.075e-02 & 1.294e-02 &  18 & 5.075e-02 & 2.808e-02 & 4.070e-02 & 5.075e-02 \\
 5 &  20 & 3.787e-02 & 1.521e-02 &  20 & 3.787e-02 & 2.783e-02 & 3.459e-02 & 3.787e-02 \\
 6 &  22 & 2.271e-02 & 1.267e-02 &  22 & 2.271e-02 & 1.943e-02 & 2.182e-02 & 2.271e-02 \\
 7 &  24 & 1.005e-02 & 6.761e-03 &  24 & 1.005e-02 & 9.157e-03 & 9.830e-03 & 1.005e-02 \\
 8 &  26 & 3.287e-03 & 2.396e-03 &  26 & 3.287e-03 & 3.069e-03 & 3.243e-03 & 3.287e-03 \\
 9 &  28 & 9.216e-04 & 6.933e-04 &  28 & 9.216e-04 & 8.749e-04 & 9.209e-04 & 9.216e-04 \\
10 &  32 & 2.639e-04 & 2.062e-04 &  32 & 2.639e-04 & 2.602e-04 & 2.631e-04 & 2.639e-04 \\
11 &  40 & 6.392e-05 & 5.842e-05 &  40 & 6.392e-05 & 6.300e-05 & 6.386e-05 & 6.392e-05 \\
12 &  54 & 9.377e-06 & 7.568e-06 &  56 & 7.955e-06 & 7.820e-06 & 7.943e-06 & 7.955e-06 \\
13 &  68 & 1.809e-06 & 1.502e-06 &  70 & 1.234e-06 & 1.222e-06 & 1.234e-06 & 1.234e-06 \\
14 &  82 & 3.144e-07 & 2.550e-07 &  84 & 1.644e-07 & 1.620e-07 & 1.641e-07 & 1.644e-07 \\
15 &  90 & 8.887e-08 & 7.358e-08 & 100 & 2.075e-08 & 2.036e-08 & 2.069e-08 & 2.075e-08 \\
16 & 102 & 1.712e-08 & 1.530e-08 & 116 & 3.001e-09 & 2.921e-09 & 2.986e-09 & 3.001e-09 \\
17 & 118 & 2.421e-09 & 1.952e-09 & 132 & 4.695e-10 & 4.549e-10 & 4.666e-10 & 4.695e-10 \\
18 & 132 & 4.695e-10 & 3.900e-10 & 144 & 9.720e-11 & 9.405e-11 & 9.715e-11 & 9.720e-11 \\
\hline
\end{tabular}
\caption{Mesh efficiency.}
\label{TabMeshEfficiency}
\end{table}

\begin{figure}
\includegraphics[width=0.7\textwidth]{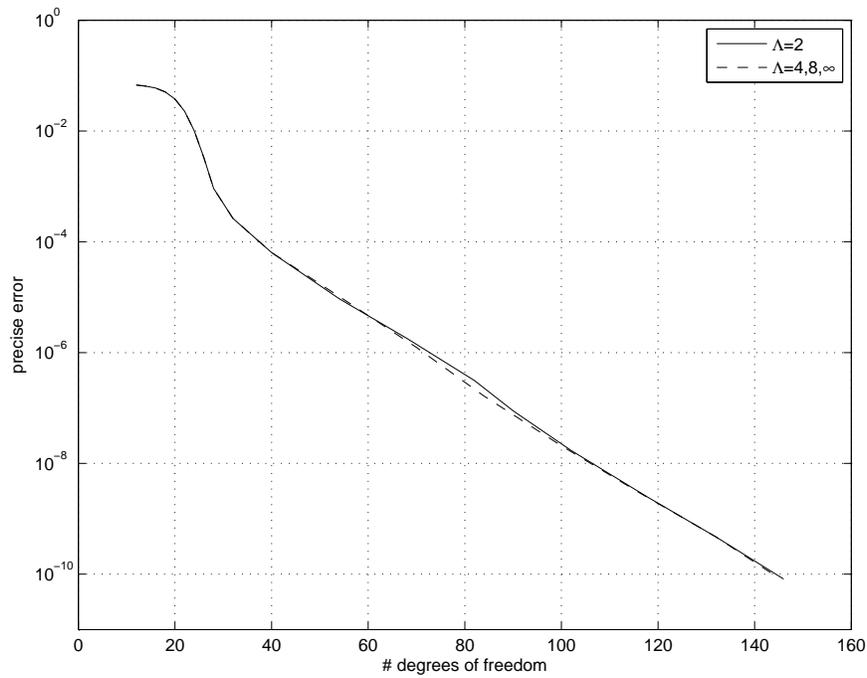}
\caption{Mesh efficiency.}
\label{FigMeshEfficiency}
\end{figure}

Table~\ref{TabMeshEfficiency} shows the number of degrees of freedom and the
corresponding error during the automatic mesh adaption for different values of
$\Lambda$. The choices $\Lambda=4$, $\Lambda=8$, and $\Lambda=\infty$ lead to
different values of $\eta$, but all result in the same mesh. For this reason,
these values share a common column for \#dof and the precise error in the table.
Figure~\ref{FigMeshEfficiency} visualizes the relationship between the number of
degrees of freedom and the error. We can see that $\Lambda=4,8,\infty$ leads to
slightly better mesh that $\Lambda=2$. However the difference is quite small so
that the faster computation with $\Lambda=2$ is very well acceptable.

\bibliographystyle{hsiam}
\bibliography{../Literature/marcel}

\end{document}